\documentclass[12pt]{amsart}

\usepackage{cancel}
\usepackage{latexsym, amssymb, amsmath}
\usepackage{soul,esint}
\usepackage{amsfonts, graphicx}
\usepackage{graphicx,color}
\usepackage{mathtools}
\usepackage{mathtools}
\usepackage{hyperref}
\usepackage{verbatim}
\makeatletter
\renewcommand*{\eqref}[1]{%
  \hyperref[{#1}]{\textup{\tagform@{\ref*{#1}}}}%
}
\makeatother

\newcommand{\be}{\begin{equation}}
\newcommand{\ee}{\end{equation}}
\newcommand{\beq}{\begin{eqnarray}}
\newcommand{\eeq}{\end{eqnarray}}

\newtheorem{thm}{Theorem}[section]
\newtheorem{conj}{Conjecture}[section]

\newtheorem{prop}{Proposition}[section]
\newtheorem{cor}{Corollary}[section]
\newtheorem{defn}{Definition}[section]
\theoremstyle{remark}
\newtheorem{rem}{Remark}[section]
\numberwithin{equation}{section}

\def\be{\begin{equation}}
\def\ee{\end{equation}}
\def\bee{\begin{equation*}}
\def\eee{\end{equation*}}

\def\la{\langle}
\def\ra{\rangle}
\def\p{\partial}

\def\tr{\operatorname{tr}}

\def\R{\mathbb{R}}

\def\S{\Sigma}
\def\SS{\partial\Omega}

\def\Lp{\Delta}
\def\Na{\nabla}
\def\eps{\varepsilon}

\def\MW{\mathcal{W}}
\def\O{\Omega}

\def\ti{\tilde}

\def\M{\mathcal{M}}
\def\mB{\mathcal{B}}
\def\mD{\mathcal{D}}

\def\mE{\mathcal{E}}

\def\i{\iota}

\begin{document}
\title[]
{Remarks on quasilocal mass and fill-ins}

 \author{Tin-Yau Tsang}
\address [Tin-Yau Tsang] {NYU Courant}
\email{tytsang@cims.nyu.edu}

\renewcommand{\subjclassname}{
  \textup{2010} Mathematics Subject Classification}

\date{Apr, 2024}

\begin{abstract} 
In this paper we would have a brief overview of several proposals of quasilocal mass which are based on Hamiltonian formulation.  We also show the positivity of the Wang-Yau energy under a more general condition.  We then further study the quasilocal mass $\mathcal{W}$ and DEC fill-ins defined in \cite{T2} in terms of completeness and shields.  
\end{abstract}


\maketitle

\markboth{TIN-YAU TSANG}{Remarks on quasilocal mass and fill-ins}
\section{Introduction}

The positive mass theorem states that an asymptotically flat initial data set $(M,g,k)$ satisfying the dominant energy condition must have non-negative ADM mass.  For the Riemannian (time-symmetric) case, it was proved by Schoen and Yau in \cite{SY1} by considering area minimisers. It is later proved by Witten in \cite{W} assuming the manifold is spin. Moreover, with \cite{PT} it is known that the Dirac operator approach works for the general case.  Schoen and Yau in \cite{SY2} proved the positive energy theorem through the Jang equation (\cite{Jang}). The positive mass theorem without the spin 
assumption was proved by Eichmair, Huang, Lee and Schoen in \cite{EHLS} by trapped surfaces. 

\

On the other hand, if a physical system is not isolated or cannot be viewed from infinity where asymptotic symmetry exists, e.g. compact manifolds with boundary, the ADM mass is not well defined. Different notions of quasilocal mass have been inspired by the Hamiltonian formulation. Those include the Brown-York mass (\cite{BY}), the Kijowski-Liu-Yau mass (\cite{KI},\cite{LY1,LY2}) and the Wang-Yau mass (\cite{WY2}).  These are all about energy perceived by a timelike observer. On the other hand, for a null observer, \cite{T2} has suggested $\mathcal{W}$ which does not need the Jang graph to show positivity. All these masses are shown to be positive under suitable energy and geometry conditions with the use of Bartnik-Shi-Tam extension \cite{B, ST} and some variants of positive mass theorem with corners along a hypersurface (\cite{M1, ST, T2}). 

\

While for the Wang-Yau mass, one of the motivations is to remove the gauge dependency of a spacelike surface. However, due to the potential blow-up of the solution to the Jang equation, some admissibility condition imposed may not be able to showcase the independence from the choice of gauge.  

\

To handle the blow-up on the Jang graph, we will show the following theorem with the tools in \cite{CZ}. 
\begin{thm}\label{RelaxedWY}(cf. \cite{WY2} Theorem 5.1)
Let $n\geq 3$. Let $(\Omega^n,g)$ be a \textbf{complete} spin Riemannian manifold with a connected boundary $\Sigma$ and suppose there exists a vector field $X$ on $\Omega$ such that 
\begin{equation} R_{g}+2div_{g} X-2|X|_{g}^2\geq 0 
\end{equation} in $\Omega$ 
and \begin{equation} H>\langle X, \nu\rangle\end{equation} on $\Sigma$, where $\nu$ is outward normal of $\Sigma$ and $H$ is the
mean curvature of $\Sigma$ with respect to $\nu$. Suppose that $\Sigma$ can be isometrically embedded into $\R^n$ as a strictly convex hypersurface. Let $H_0$ be the mean curvature of the isometric embedding of $\Sigma$ into $\R^n$ with respect to the normal pointing to the infinity. Then
\be \label{maininequality}
\int_\Sigma H_0  \geq \int_\Sigma H-\langle X, \nu\rangle.
\ee
\end{thm}

As a result, we can relax the admissibility condition imposed in \cite{WY2} to show the Wang-Yau energy is non-negative on a bigger set of isometric embeddings of surfaces into Minkowski space. Also note that, we at the same time prove a generalisation of the time-symmetric charged positive mass theorem (cf. \cite{GHHP}) and a weighted positive mass theorem (\cite{BO, CZ}) with arbitrary complete ends. 

\

For geometry on compact manifolds with boundary, Gromov proposed the following conjecture. 
\begin{conj}\label{Gromov}(\cite{Gro19} Sect 3.12.2 III., IV.) Let $(M, g)$ be a compact Riemannian manifold with scalar curvature $R \geq \sigma$. Then there exists $\Lambda$ depending only on $\sigma$ and the intrinsic geometry of $(\p M, g|_{T(\p M)})$ such that
\be
\int_{\p M} H \leq \Lambda,
\ee
where $H$ is the mean curvature of the boundary $\p M$ in $(M, g)$ with respect to the outward unit normal vector.  
\end{conj}

In \cite{SWWZ} and \cite{SWW}, there was a partially affirmative answer given by the parabolic method to construct a metric extension done in \cite{ST}. While in \cite{T2}, a corresponding conjecture regarding energy for initial data sets are proposed and partially confirmed in terms of the Bartnik data sets (\cite{B2}) and the idea of fill-ins satisfying dominant energy condition.  

\ 

Along with the recent development of positive mass theorem of arbitrary ends (\cite{LUY2,LLU2,Zhu,CLSZ,CZ}), \cite{LLU3} proved a Riemannian positive mass theorem with corners with arbitrary complete ends and shields.  It is then used to show the positivity of the Brown-York mass of the boundary of manifolds which are complete or with shields.  While \cite{CLZ} adopted the Dirac operator approach to show a positive mass theorem with corners with arbitrary complete ends and shields for initial data sets.  Correspondingly, $\MW$ and DEC fill-ins defined in \cite{T2} can be generalised to take those into considerations as well.  

\

This text is organised as follows. In Section \ref{Preliminaries}, we recall asymptotically flat initial data sets and the Hamiltonian formulation.  In Section \ref{review}, several proposals of quasilocal mass are reviewed. In Section \ref{PWY}, we first review the admissiblility conditions in \cite{WY2}. Then,  we will discuss the strategy of showing that the Wang-Yau mass is non-negative and how the admissibility requirement can be relaxed.  The positivity of $\MW$ and the non-existence of DEC fill-ins with completeness and shields is discussed in Section \ref{WDECfillin}.  

\

\textbf{Acknowledgements }: The author would like to thank Prof.  Mu-Tao Wang for some inspiring discussion on \cite{WY0}.   The author is grateful to Prof. Ricahrd Schoen and Kai-Wei Zhao for explaining \cite{SY2} in detail. The author thanks Prof. Lan-Hsuan Huang and Prof. Dan Lee for sharing their insights on the improvability of energy.  Moreover, the author would like to thank Prof. Chao Li, Chi-Cheuk Tsang and Douglas Stryker for helpful discussion.   

Finally, the author is grateful to the National Centre for Theoretical Sciences (NCTS), in particular Prof. Yng-Ing Lee for her generous support, for organising the 5th Taiwan International Conference on Geometry in honour of Prof. Schoen's 70th birthday where part of this manuscript is done and discussed.  

\section{Preliminaries}\label{Preliminaries}
\subsection{Asymptotically flat initial data sets}
Given an initial data set \\ $(M^n,g,k)$, where $g$ is a Riemannian metric and $k$ is a symmetric $(0,2)$-tensor. Define the conjugate momentum tensor by $\pi=k-(tr_g k)g$. Under constraint equations, we can define the mass density $\mu$ and the current density $J$ by
$$\mu=\frac{1}{2}(R_g+(tr_g k)^2-|k|_g^2)$$
and $$J=div_g(k-(tr_gk)g)=div_g\pi.$$
$(M,g,k)$ is said to satisfy the dominant energy condition if 
$$\mu\geq|J|_g.$$ 
We say $(M, g, k)$ is asymptotically flat if there exists a compact set $\mathcal{C}\subset M$ such that  $M\setminus \mathcal{C}=\coprod_{i=1}^k N_i$, where each end $N_i=\R^3\setminus B_{r_i}$ through a coordinate diffeomorphism in which 
$$
g_{ij} = \delta_{ij} + O^2(|x|^{-q}),
$$
and 
$$
k_{ij}=O^1(|x|^{-q-1}),
$$
where $q>\frac{n-2}{2}$, $ \mu, J \in L^1(M)$ and for a function $f$ on $M$, $f = O^{m}(|x|^{-p})$ means $\sum_{|l|=0}^m||x|^{p+|l|}\partial^lf|$ is bounded near the infinity.

\

For each end, the ADM energy-momentum vector $(E,P)$ and the ADM mass $\mathfrak{m}$ \cite{ADM} are given by 
$$
E = \frac{1}{2(n-1)\omega_{n-1}} \lim_{r \to \infty} \int_{|x| = r}  ( g_{ij,i} - g _{ii, j} ) \nu^j,
$$ 
$$P_i:=\frac{1}{(n-1)\omega_{n-1}} \lim_{r \to \infty} \int_{|x| = r}  \pi_{ij} \nu^j, \,\ \,\ \,\ i=1,2,...,n$$ 
and 
$$\mathfrak{m}=\sqrt{E^2-|P|^2},$$
where the outward unit normal $\nu$ and surface integral are with respect to the Euclidean metric and $\omega_{n-1}=|\mathbb{S}^{n-1}|_{g_{Euc}|_{S^{n-1}}}$.  The well-definedness is shown in \cite{Bartnik}. 

For a smooth closed hypersurface $\S\subset M$, we say $\S$ is a weakly outer trapped surface if on $\S$, the outer null expansion 
$$\theta_+=H+tr_\S k \leq 0,$$
and a marginally outer trapped surface ($MOTS$) if 
$$\theta_+=0;$$
correspondingly, 
a weakly inner trapped surface if the inner null expansion 
$$\theta_-=H-tr_\S k \leq 0,$$
and a marginally inner trapped surface ($MITS$) if 
$$\theta_-=0,$$ 
where $H$ is computed with respect to the normal pointing to the infinity of the designated end $\mathcal{E}$. A surface is weakly trapped if it is either weakly outer trapped or weakly inner trapped. 
In this note, the mean curvature is computed in the convention that $\mathbb{S}^2\subset \R^3$ has positive mean curvature 2 with respect to the outward normal.  

For notations, in the case $k\equiv 0$, we denote $(M,g,0)$ by $(M,g)$ and it is called an asymptotically flat manifold. For a spacelike $(n-1)$-surface $\S$ in an $(n+1)$-spacetime $N$, the mean curvature vector is denoted by $\vec{H}$. Let $\{e_n,e_{n+1}\}$ be an orthonormal frame of $N\S=TN/T\S$, where $e_n$ is spacelike and $e_{n+1}$ timelike.  Let $k$ and $A$ be the second fundamental form of $\S$ with respect to  $e_{n+1}$ and $e_n$ respectively.  And $H:=tr_{\S} A$.  Therefore, $\vec{H}:=tr_{\S} k \, e_{n+1}-H e_n$. 

\subsection{Hamiltonian formulation (Hamilton-Jacobi analysis)}\label{Hamiltonian formulation}
Let $(\Omega^n,g,k)$ be a compact initial data set with boundary $\S$.  A spacetime $(N^{n+1},\bar g)$ with boundary $\bar \S$ can be constructed by infinitesimally deforming the initial data set $(\Omega,g,k,\S)$ in a transversal, timelike direction $\p_t= V \vec{n} + W^i\p_i$ which satisfies $\bar\Na_{\p_t} t=1$, where $V$ is the lapse function, $\vec{n}$ is the future pointing timelike unit normal of $\Omega$ in $N$ and $W$ is the shift vector. Further assume that $\Omega$ meets $\bar \S$ orthogonally. The purely gravitational contribution $\mathcal{H}_{grav}$ to the total Hamiltonian at the slice $\Omega$ is given by (\cite{ADM},\cite{RT},\cite{BY},\cite{HH})
\be \label{Ham}
c(n) \mathcal{H}_{grav}(V,W) = \int_{\Omega} (\mu V  + \la J,W \ra) - \int_{\S}(HV-\pi(\nu,W)),
\ee
where $H$ is the mean curvature of $\S$ with respect to the outward normal $\nu$ of $\Omega$ and $\pi$ is the conjugate momentum tensor.  From this, the mean curvature $H$ and the 1-form $\pi(\nu,\cdot)$ can be seen as energy perceived by an observer on the boundary. 

\section{Review on proposals of quasilocal mass based on Hamiltonian formulation}\label{review}
In this section, we would introduce several proposals of quasilocal mass whose principle is the same, Hamiltonian formulation. From \eqref{Ham} in Section \ref{Hamiltonian formulation}, if we are only given $(\Omega,g,k)$,  the boundary integrand $HV-\pi(\nu,W)$ can be considered as $\bar g ( \vec{E}_{\S,\nu} , \p_t)$, where $\vec{E}_{\S,\nu}:= -H \vec{n} +\pi(\nu, \cdot)$, the boundary energy vector. There are freedom for us to choose the time like observer $\p_t$, i.e $V$ and $W$.  On the other hand, if we see $\S$ as a spacelike codimension 2 surface in a spacetime $(N,\bar g)$, there is extra freedom for us to pick a frame of the 2 dimensional normal bundle $N\S=TN/T\S$.  The boundary energy vector $\vec{E}_{\S,\cdot}$ would change since the second fundamental forms $k$ and $A$ would change according to the timelike normal direction and the spacelike normal direction respectively. (From this perspective, giving $\Omega$ is the same as choosing a gauge by fixing $e_n=\nu$.) These freedoms lead to different suggestions of quasilocal mass.  Nonetheless, they are all defined to be the difference of the boundary integral in \eqref{Ham} in a given initial data set (in a spacetime) and that of a reference.  

\subsection{The Brown-York mass \cite{BY} }Let $(\Omega, g)$ be a compact manifold with boundary $\S$. Let $H$ denote the mean curvature with respect to the outward normal. If $\S$ can be isometrically embedded into the Euclidean space, then let $H_0$ be the mean curvature of $\S$ with respect to the normal pointing to $\infty$ in $\R^n$. 
\be
m_{BY}:=\int_{\S} H_0-H . 
\ee
In this case, $V=1$ and $W=0$, the reference is $\R^n$ while the gauge is given by $\Omega$. 
\subsection{The Kijowski-Liu-Yau mass} As a generalisation of the Brown-York mass, Kijowski (\cite{KI}) and Liu and Yau(\cite{LY1,LY2}) independenty proposed the following. Let $(\Omega, g,k)$ be a compact initial data set with a boundary $\S$. If $\S$ can be isometrically embedded into the Euclidean space, then let $H_0$ be the mean curvature of $\S$ with respect to the outward normal in $\R^n$,  assuming $H>|tr_{\S}k|$. 
\be
m_{KLY}:=\int_{\S} H_0-|\vec {H}|, 
\ee
i.e.  the reference is $\R^n(\subset \R^{n,1})$ while for the gauge and the choice of observer, we shall first see the following discussion of Wang-Yau mass.  
\subsection{The Wang-Yau mass} In \cite{WY0,WY2}, in order to take the reference as a spacelike slice of Minkowski space, a gauge dependent generalised mean curvature is first defined. Let $\iota:\S \hookrightarrow (N,\bar g)$. For each frame $\{e_n, e_{n+1}\}$ of $N\S$ where $e_n$ is spacelike and $e_{n+1}$ is future pointing timelike,  let $A$ and $k$ to denote the second fundamental form with respect to $e_n$ and $e_{n+1}$ respectively.  For any $\tau:\S\to \R$, define 
\be
\begin{split}
h(\iota, \tau, e_n):=&(\sqrt{1+|\Na^{\S} \tau|^2})H+ k(e_n, \Na^{\S} \tau)\\
=&-(\sqrt{1+|\Na^{\S} \tau|^2})\bar g (  \vec H, e_n )+ \pi(e_n, \Na^{\S} \tau)\\
=& \bar g (\vec{E}_{\S,e_n},\sqrt{1+|\Na^{\S} \tau|^2}e_{n+1}+\Na^{\S}\tau ).
\end{split}
\ee

Thus, $V= \sqrt{1+|\Na^{\S} \tau|^2}$ and $W=\Na^{\S} \tau$.  

Then, assuming $H>|tr_{\S} k|$, one can define a quantity 
\be 
\mathfrak{H}(\iota, \tau):=\int_{\S} \inf_{e_n} h(\iota, \tau, e_n). 
\ee 

Further assume there exists an isometric embedding $\i_\tau^0: (\S, \i^*\bar g) \hookrightarrow \R^{n,1}$ with $\tau$ being the height with respect to a constant time slice,  the Wang-Yau energy is defined as follows. 
\be 
E_{WY}(\tau):=\mathfrak{H}(\iota_{\tau}^0, \tau)-\mathfrak{H}(\iota, \tau).  
\ee 

Then, one can in principle define a mass by $\inf_{\tau}E_{WY}(\tau)$ while one can only show $E_{WY}(\tau)\geq 0$ if $\tau$ satisfies the admissibility conditions (\cite{WY2} Defintion 5.1). which are to be further discussed in the next section.  Let $\mathcal{C}$ be the collection of admissible functions on $\S$.  Finally, the Wang-Yau mass is defined as
\be 
m_{WY}:=\inf_{\tau\in\mathcal{C}}E_{WY}(\tau). 
\ee 

Note that the Kijowski-Liu-Yau mass $m_{KLY}$ is taking $\tau=0$ and the gauge with $e_n=-\frac{\vec{H}}{|\vec{H}|}$ which is not necessarily $\nu\in T\Omega$. Moreover, $V=1$ and $W=0$.  Coincidentally under this gauge, the quantity considered only needs information from the given  initial data set instead of a spacetime.  
\subsection{$\MW$} 
All the previous masses are considering a unit timelike observer.  Inspired by the spacetime positive mass theorem with corners, in \cite{T2} a null observer is used to define a quasilocal mass. Let $(\Omega, g,k)$ be a compact initial data set with boundary $\S$. If $\S$ can be isometrically embedded into the Euclidean space, then let $H_0$ be the mean curvature of $\S$ with respect to the outward normal in $\R^n$,   
\be
\mathcal{W}:=\int_{\S} H_0-(H-|\pi(\nu,\cdot)|). 
\ee
Therefore, the reference is $\R^n$, the gauge is fixed by $\Omega$ while $(V,W)$ is formally $(1,1)$. 

\begin{rem} Recently, a quasilocal mass was suggested in \cite{AKY2} which also used a null observer while considering the extra freedom of gauge as in the Wang-Yau energy approach. 
\end{rem}

\section{Positivity of the Wang-Yau mass} \label{PWY}
In this section, we first discuss the idea presented in \cite{WY0,WY2}. This is combining the work of the Jang equation and the proof of positivity of the Brown-York mass in \cite{ST}. 

\subsection{Admissibility in \cite{WY2}} 
Let $n=3$ (see Remark \ref{HD} below). Let $\i:\S \hookrightarrow (N,\bar g)$ as a spacelike connected $(n-1)$-surface in a $(n+1)$-spacetime. 
In order to show the Wang-Yau mass is non-negative, one only has to show the Wang-Yau energy is non-negative for each admissible $\tau$. Admissibility conditions of $\tau$ means  
\begin{enumerate}
\item \label{1} the existence of $\i^0_{\tau}$ and $(\S, (\i^0_{\tau})^*\eta+d\tau^2)$ is a strictly convex hypersurface in $\R^n$, where $\eta$ denotes the Minkowski metric.  
\item \label{2} the existence of a compact spacelike slice $\Omega \subset N$ which satisfies the dominant energy condition such that $\S=\SS$ with $H>|tr_{\S}k|$ over which the Jang equation can be solved with boundary value $\tau$ with bounded gradient near $\S$, and 
\item \label{3} the positivity of $h(\i,\tau,e_n)$ for a specified frame $\{e_n,e_{n+1}\}$ of $N\S$ (See \cite{WY2} Theorem 4.1).    
\end{enumerate}
Item \eqref{1}'s first requirement allows the definition of $E_{WY}(\tau)$ and item \eqref{2} allows us to consider the Jang graph of $\Omega$ instead of $\Omega$. Let $\ti \Omega$ denote the Jang graph of $\Omega$. The second requirement in item \eqref{2} and item \eqref{3} allows the use of the Bartnik-Shi-Tam extension \cite{ST,B3}) which is asymptotically flat and scalar flat with prescribed mean curvature along $\p \ti \O$. At the end, by gluing $\ti \Omega$ and the extension together, one can use a variant of positive mass theorem with corners for spin manifolds to conclude that $E_{WY}(\tau)$ is non-negative by comparing with the ADM mass of the glued manifold using the monotonicity discovered in \cite{ST} Lemma 4.2.

\begin{rem}\label{HD} There are 2 ways to consider $E_{WY}$ in a higher dimension. 
\begin{enumerate}
\item If one assume item \eqref{2}, then one can use the mollification technique in \cite{AKY1} (cf. \cite{M1}) to prove the theorem for the dimension of $\S$ up to 6.  
\item If one further imposes in item \eqref{2} that the spacelike slice is spin, then one can show positivity for all dimensions. 
\end{enumerate}
\end{rem}

\begin{rem}
One can also impose the spacetime dominant energy condition on the spacetime $N$ so that every spacelike slice inherits the dominant energy condition automatically. 
\end{rem}


\

We shall have a closer look at item \eqref{2} regarding the solvability of the Jang equation over the whole $\Omega$.  From \cite{SY2} Section 4 Proposition 4, we know the solution to the Jang equation $f$ may blow up ($|f|\to\infty$). When it does, this takes place at the apparent horizon $\S_{AH}$ which consists of MOTS and MITS.  Let $\S_B\subset \S_{AH}$ denote the component of the apparent horizon where the blow up takes place. If the base initial data set $\Omega$ satisfies the strict dominant energy condition, then we know $\S_{B}$ is Yamabe positive and \cite{SY2} Section 5 is able to use the first eigenvalue of the operator $-\Lp_{\hat g} + \frac{(n-2)}{4(n-1)}R_{\hat g}$ which is strictly positive, where $\hat g=(\bar g|_{\Omega})|_{\S_B}$ to close up the asymptotic cylinder $\Sigma_{B}\times$ half-line to make it equivalent to a ball and with non-negative scalar curvature. (See some further details in \cite{E3} Section 3.) However, this method does not work if $\S_{B}$ is not Yamabe positive.  Therefore, item \eqref{2} is actually not only a restriction on the choice of function $\tau$, but also somehow a restriction on $N$ which allows the existence of suitable spacelike slices.  

\begin{rem}
There is another approach in \cite{Mon} by spinors to show the Brown-York mass is non-negative without using Bartnik-Shi-Tam extension. This has been generalised in \cite{MY} to show the positivity of Wang-Yau mass under the admissibility requirement above.  
\end{rem}

\subsection{Relaxing admissibility}
As aforementioned, item \eqref{2} could be too restrictive as it somehow also limits the spacetime and spacelike slices we can choose.  In this part, we would discuss some recent development of positive mass theorems which can relax the admissibility conditions imposed in \cite{WY2}. 

\

Following up with the the procedure of closing up the cylinder, as motivated by the proof of positive mass theorem in \cite{SY1, EHLS}, one may want to perturb the dominant energy condition to the strict dominant energy condition. Moreover, since we are considering both the intrinsic geometry and extrinsic energy of $\S$, we want to keep them the same along the change. But it is not always possible, see \cite{Cor1,CH,HL2}.  Physically speaking, one shall not get energy for free.  For example, by the rigidity result of the positive mass theorem with corners \cite{M1, ST},  it is impossible to deform a standard unit ball to have strictly positive scalar curvature while keeping the boundary isometric to a standard unit sphere with the mean curvature equal to a unit standard sphere in the Euclidean space.  Geometrically speaking,  for example, we can tell from the Gauss-Bonnet Theorem on 2 surface that the interior curvature change would easily affect the boundary geometry.  

\

In a very recent work \cite{AYZ},  assuming $\S_B$ is a strictly stable MOTS/ MITS,  the principal eigenvalue of the stability operator (the first variation of $\theta_{\pm}$) is used instead of that of $-\Lp_{\hat g} + \frac{(n-2)}{4(n-1)}R_{\hat g}$.  That eigenvalue provides a parameter to approximate the rate of convergence of the blow up ends to an exact cylinder. Under this more specific understanding of the asymptotics, one can establish the existence of harmonic spinors to allow the use of positive mass theorem with corners. 

\

Here, we provide another solution so that we do not need to care about the MOTS/ MITS stability of $\S_B$.  Thus, we can further relax item (2).  Denote the region in $\Omega$ outside of $\S_B$ by $\O_B$, i.e. $\p\Omega_B=\Sigma\cup\S_B$.  It is noticed that by the Harnack inequality in Proposition 2 in \cite{SY2}, we know that the Jang graph $\ti\Omega_{B}:=f(\O_B)$ is complete.  Note that as shown in \cite{WY2} Theorem 4.2, $H>|tr_{\S} k|$ already ensures the gradient of the solution would not blow up at $\S$.   Denote the gluing of $(\ti \O_B, \ti g:=(\bar g|_{\O_B}+dt^2)|_{\ti \O_B})$ and the Bartnik-Shi-Tam extension $(M_{BST},g_{BST})$ along $\S$ by $(M,g)$.  $(M,g)$ possesses an asymptotically flat end obviously, and as discussed, with other asymptotically cylindrical ends ($\S_B\times$ half-line) under the presence of $\S_B$.  $M$ is thus not an asymptotically flat manifold. 

On the other hand,  positive mass theorem with arbitrary complete ends is proved in \cite{LUY2, CZ, LLU2, Zhu, CLSZ}.  In short, it sates that any asymptotically flat end of a complete manifold with non-negative scalar curvature must have non negative ADM mass.  Note that $M$ is a 3 manifold which is always spin while the scalar curvature on $\ti\O_B$ satisfies the inequality (\cite{SY2} (2.29)) for some vector field $X$, 
\be
R_{\ti g}+2 div_{\ti g} X -2|X|_{\ti g}^2 \geq 0. 
\ee
Therefore, we choose to follow the spinorial approach in \cite{CZ} to generalise the proof of Theorem 5.1 in \cite{WY2} under the presence of blow up.  In particular, we show Theorem \ref{RelaxedWY}.  For reading convenience,  we restate the theorem here.   

\begin{thm}
Let $n\geq 3$. Let $(\Omega^n,g_-)$ be a complete spin Riemannian manifold with a connected boundary $\Sigma$ and suppose there exists a vector field $X$ on $\Omega$ such that 
\begin{equation}\label{jang_scalar} R_{g_-}+2div_{g_-} X-2|X|_{g_-}^2\geq 0 
\end{equation} in $\Omega$ 
and \begin{equation}\label{jang_mean} H_->\langle X, \nu\rangle\end{equation} on $\Sigma$, where $\nu$ is outward normal of $\Sigma$ and $H_-$ is the
mean curvature of $\Sigma$ with respect to $\nu$. Suppose that $\Sigma$ can be isometrically embedded into $\R^n$ as a strictly convex hypersurface. Let $H_0$ be the mean curvature of the isometric embedding of $\Sigma$ into $\R^n$ with respect to the normal pointing to the infinity. Then
\be \label{maininequality}
\int_\Sigma H_0  \geq \int_\Sigma H_--\langle X, \nu\rangle.
\ee
Moreover, the equality is achieved if $(\Omega,g_-)$ is isometric to a domain in $\R^n$ and $X=0$. 
\end{thm}

\begin{proof}
We notate the image of isometric embedding of $\S$ into $\R^3$ by $\S_{0}$, and the unbounded region of $\R^3$ outside of $\S_{0}$ by $M_+=\S_{0} \times [0,\infty)$ which stands for a foliation by unit normal flow.  
Then as in \cite{ST}, we can construct an asymptotically flat metric $g_{+}=u(r)^2dr^2+g_r$ with zero scalar curvature on $M_{+}$ (\cite{ST} Theorem 2.1(b)), where $u(0)=u>0$ and $g_r$ stands for the metric induced on $\S_r=\Sigma_{0} \times \{t=r\}$ by the Euclidean metric on $\R^n$. Since $\S_r$ is strictly convex, we have (Lemma 4.2 in \cite{ST}), 
\be
\begin{split}
(n-1)\omega_{n-1} \frac{d}{dr}Q(\S_r)
:=&\frac{d}{dr} \int_{\S_r}H_0(r)\left( 1-\frac{1}{u(r)} \right)d\sigma_r \\
=&-\frac{1}{2} \int_{\S_r} R_{\S_r} u^{-1}(1-u)^2\leq 0, 
\end{split}
\ee
where $H_0(r)$ is the mean curvature of $\S_r$ with respect to the Euclidean metric of $\R^3$. Moreover, by Theorem 2.1 (c) in \cite{ST}, we have$$\lim_{r\to \infty}Q(\S_r)= E(g_{+}).$$ Therefore to show \eqref{maininequality}, it suffices to show $E(g_{+})\geq 0$ under a choice of $u$ such that $$Q(\S_0)=\frac{1}{(n-1)\omega_{n-1}} \int_{\S}H_0-(H_--\la X,\nu \ra),$$ i.e.

\be
u=\frac{H_0}{H_--<X,\nu>}.
\ee 

Moreover, under this choice of $u$, let $H_+$ denote the mean curvature of $\S_0$ with respect to the normal pointing to $\infty$ of $M_+$, then $H_+=H_--\la X,\nu \ra$.  

Let $(M,g):=(\Omega \cup M_+, g_-\cup g_+)$. For $j\in \mathbb{N}$, let $\S\subset \Omega_j \subset \Omega_{j+1}\subset \Omega$,  such that $\{\Omega_j\}_{j>0}$ is a compact exhaustion of $\Omega$ with smooth boundary.  Denote $\p\Omega_j \setminus \S$ by $\S_j$.  Let $(M_j,g_j):=(\Omega_j \cup M_+, g_-|_{\Omega_j}\cup g_+)$.

\

We would use the set up of the spinor bundle $S$,  Callias operator $\sigma$ and Penrose operator $\mathcal{P}$, the chirality $\chi$ on boundary $\chi u:=\nu \cdot \sigma u$ and the definition of asymptotically constant spinor introduced in \cite{CZ} Section 2. Without loss of generality, by scaling, we can assume the constant spinor to which a spinor section being asymptotic is with norm 1 . 

Let $\psi$ be a compactly supported Lipschitz function on $M$.  Let $\mB_{\psi}:=\mathcal{D}+\psi \sigma$, where $\mD$ is the Dirac operator. Define $\theta_{\psi}^X$, $\eta_{\psi}^X$, $\bar \theta_{\psi}^X$ and $\bar \eta_{\psi}^X$ as follows. 
\begin{enumerate}
\item On $\Omega$, 
\begin{enumerate}
\item $\theta_{\psi}^X=\frac{1}{4}\left( R_g+2div_g X- 2|X|_g^2 \right)+\psi^2-|\Na \psi|$.
\item $\bar \theta_{\psi}^X=\frac{n}{2(2n-1)}\left( R_g+2div_g X- 2|X|_g^2 \right)+\psi^2-|\Na \psi|$.
\end{enumerate}
\item On $M_+$, 
\begin{enumerate}
\item $\theta_{\psi}^X=\frac{1}{4}R_g+\psi^2-|\Na \psi|$.
\item $\bar \theta_{\psi}^X=\frac{n}{2(2n-1)} R_g+\psi^2-|\Na \psi|$.
\end{enumerate}
\item On $\S_j$
\begin{enumerate}
\item $\eta_{\psi}^X=\frac{1}{2}(H+\la X,\nu\ra)+\psi$.
\item $\bar \eta_{\psi}^X=\frac{n}{2n-1}(H+\la X,\nu\ra)+\psi$,
\end{enumerate}
where $H$ denotes the mean curvature with respect to the outward pointing normal $\nu$. 
\end{enumerate}

We now proceed the following.  
\begin{prop}\label{spinor}(cf. \cite{CZ} Proposition 2.5, \cite{WY2} (5.5) )
Let $u\in W^{1,2}_{loc}(M_j, S)$ with $u=\chi u$ on $\S_j$.  Further assume that $u$ is asymptotically constant.  Then we have, 
\be \label{crucial1}
\frac{n-1}{2}\omega_{n-1}E(g_+)+||\mB_{\psi} u||_{L^2(M_j)}^2 \geq \frac{1}{2}||\Na u||_{L^2(M_j)}^2+\int_{M_j} \theta_{\psi}^X |u|^2+\int_{\S_j} \eta_\psi^X |u|^2. 
\ee
and
\be \label{crucial2}
\begin{split}
\frac{n(n-1)}{2n-1}\omega_{n-1}E(g_+)+||\mB_{\psi} u||_{L^2(M_j)}^2\geq \frac{n}{2n-1}||\mathcal{P} u||_{L^2(M_j)}^2+\int_{M_j} \theta_{\psi}^X |u|^2+\int_{\S_j} \eta_\psi^X |u|^2. 
\end{split}
\ee
\end{prop}
\begin{proof}
First by Cauchy-Schawrz inequality,  
\be \label{CauchySchwarz}
X(|u|^2)=\langle \nabla_Xu, u\rangle
+\langle u, \nabla_X u\rangle\geq
-2|\nabla_X u||u|\geq -|\Na u|^2-|X|^2|u|^2. 
\ee
Therefore, 
\be \label{scalaruseful}
\frac{1}{2}|\nabla u|^2+\frac{1}{2}X(|u|^2)\geq -\frac{1}{2}|X|^2|u|^2.
\ee

Integrating by parts on $\Omega_j$, we obtain
\be \label{bypart}
\int_{\S_j}\langle X, \nu\rangle|u|^2+\int_{\S}\langle X, \nu\rangle
|u|^2=\int_{\Omega_j} div_{g} X|u|^2+\int_{\Omega_j} X(|u|^2).
\ee
Then,  one can integrate the Lichnerowicz formula (\cite{Lichnerowicz}) on $\Omega_j$ and $M_+$ respectively and sum them up. Then use the fact that $H_+=H_--\la X,\nu \ra$ on $\S$ with \eqref{scalaruseful} and \eqref{bypart} to obtain the result. 
\end{proof}
By Proposition \ref{spinor}, one can establish the invertibility of $\mB_{\psi}$ as in \cite{CZ} Section 2. And then, one can proceeds the proof of $E(g_+)$ is non negative as in \cite{CZ} Section 3 and 4 by considering the exhaustion $(M_j,g_j)$ as $j\to \infty$.  Rigidity comes from the existence of parallel spinors (as we have freedom to choose the constant spinor reference) and $X=0$ by \eqref{CauchySchwarz}.

\end{proof}

\begin{cor}\label{MAIN}
For $n=3$, correspondingly, we can relax item (2) in the admissibility conditions to allow the blow up of the solution to the Jang equation in $int \, \O$.  
\end{cor} 

\begin{rem} Note that by naturality, any open subset of a spin manifold is still spin. If $\Omega$ is further assumed to be spin in item (2), and if there are blow-ups of the solution to the Jang equation in $\Omega$, we can only handle $E_{WY}$ with $3 \leq n\leq 7$.  The obstruction comes from the regularity of $\S_B$ and thus application of a suitable positive mass theorem with corners.  
\end{rem}

\begin{rem} Let $(M,g)$ be a complete spin manifold and $\mE$ be an asymptotically flat end of $M$. From the proof of Theorem \ref{RelaxedWY}, one can see that we have generalised and/ or provided another perspective and proof of the following. 
\begin{enumerate}
\item If $X$ is a asymptotically vanishing divergence free electric field, then one can prove a time-symmetric charged positive mass theorem on $\mE$ if the charged dominant energy condition is satisfied (c.f. \cite{GHHP}).  
\item If $X=\Na f$, where $f$ is asymptotically vanishing, then one can show a weighted positive mass theorem with the use of one half of \eqref{scalaruseful} if Perelman's scalar curvature is non-negative (c.f. \cite{BO, CZ23}). Moreover,  our conclusion of rigidity is stronger than the one shown in \cite{BO}.  
\end{enumerate} 
\end{rem}

\begin{rem}
The corresponding results for shield type statements (see Section \ref{WDECfillin} Definition \ref{dec.shield}) is achieved exactly in the same way. 
\end{rem}

\section{$\mathcal{W}$ and DEC fill-ins with completeness and shields }\label{WDECfillin}
In \cite{T2}, it was proved that $\mathcal W$ is non-negative on $\S$ given that $\S=\SS$, where $(\Omega,g,k)$ is a spin compact initial data set satisfying the dominant energy condition with $H>|\pi(\nu,\cdot)|$ on $\S$ and $\S$ can be isometrically embedded into $\R^n$ as a strictly convex hypersurface, by the spacetime positive mass theorem with corners on spin initial data sets which can be achieved by summing up the integration of Lichnerowicz formula on $\Omega$ and its corresponding extension $M_{BST}$.  (See \cite{ST}, \cite{HL1} Section 7, also \cite{Shibuya}, \cite{LL}. ) 

Note that the above is done without using the Jang graph. Therefore, it can be shown to be non-negative for spin initial data sets without the dimension restriction arises from the regularity of $\S_B$ for the  blow-up of the solution to the Jang equation.  

Moreover, in \cite{T2}, it is shown that by spacetime positive mass theorem with corners, if $\S$ is having too much energy (in various ways), then it cannot admit a compact spin fill-in satisfying the dominant energy condition.  Now let us recall the definition of a dominant energy shield in \cite{CLZ}.  
\begin{defn}\label{dec.shield}
Let $(M^n,g,k)$ be an initial data set, not assumed to be complete. We say that $(M,g,k)$ \emph{contains a dominant energy shield $U_0 \supset U_1 \supset U_2$} if $U_0$, $U_1$, and $U_2$ are open subsets of $M$ such that 
$U_0 \supset \overline U_1$, $U_1\supset \overline U_2$,
 the closure of $U_0$ in $(M,g)$ is a complete manifold with compact boundary,
and we have the following: 
\begin{enumerate}
    \item $\mu-|J|\ge 0$ on $U_0$,
    \item $\mu-|J| \ge \sigma n(n-1)$ on $U_1\setminus U_2$,
    \item the mean curvature $H_{\partial\bar U_0}$ on $\partial\bar U_0$ and the symmetric two tensor $k$ satisfy
    \[H_{\partial\bar U_0} -\left|\pi(\nu,\cdot)|_{T{\p\bar U_0}}\right|>-\Psi(d,l).\]
    Here, $\Psi(d,l)$ is the constant defined as
    \[
    \Psi(d,l) := \begin{cases}  
    \frac{2}{n} \frac{\lambda(d)}{1-l \lambda(d)} & \text{if \(d < \frac{\pi}{\sqrt{\sigma} n}\) and \(l < \frac{1}{\lambda(d)}\),} \\
    \infty & \text{otherwise,}
    \end{cases} 
    \]
    where $d := dist_g(\partial U_2, \partial U_1)$, $l := dist_g(\partial U_1, \partial U_0)$, and
    \[
    \lambda(d):= \frac{\sqrt{\sigma} n}{2} \tan\left(\frac{ \sqrt{\sigma} n d}{2}\right).
    \]
\end{enumerate}
\end{defn}

Correspondingly, we an extend the definition of fill-ins as follows.  
\begin{defn}(cf. \cite{B2} Definition 2, \cite{SWWZ})
For $n\geq 3$, a tuple $(\Sigma^{n-1}, \gamma, \alpha, H, \beta)$ is called a spacetime Bartnik data set, where $(\Sigma, \gamma, \alpha)$ is an oriented closed null-cobordant initial data set with $\alpha\in C^{1,\eps}$, while $H$ and $\beta$ are respectively a smooth function and a $C^{1,\eps}$ 1-form on $\Sigma$. A compact or \textbf{complete} initial data set $(\Omega^n,g,k)$ is called a fill-in of $D_{SB}$ if there is an isometry $\phi: (\Sigma^{n-1}, \gamma)\to (\p\Omega, g|_{\p\Omega}) $ such that 
\begin{enumerate}
\item $\phi^*H_g=H$,  where $H_g$ is the mean curvature of $\p\Omega$ to $g$ with respect to the outward unit normal $\nu$,  
\item $\phi^*tr_{\p \Omega}k=tr_{\S}\alpha$,  and
\item $\phi^*( k(\nu,\cdot))=\beta$.
\end{enumerate} 
Correspondingly,  we can also consider the case where an initial data set $(\Omega,g,k)$ is with a domain energy shield, if there exists an isometry $\phi:(\Sigma^{n-1}, \gamma)\to (\p\Omega_{\S}, g|_{\p\Omega_{\S}})$ such that above condition holds with $\p\O_{\S}\subset U_2$ and $\p\Omega\setminus \SS_{\S}= \p\bar U_0 \sqcup \S_{AH}$, where $\S_{AH}$ only consists of MOTS/MITS. 
\end{defn}
We can see from the above definition that on $\S$,  $$\phi^*(|\pi(\nu,\cdot)|_g) =\sqrt{(\tr_{\S}\alpha)^2+|\beta|_{\gamma}^2}.$$ 

On the other hand,  \cite{LLU3} has shown that $m_{BY}$ is non-negative on surfaces $\S$ which is a boundary of a complete manifold $(\Omega,g)$ with $R_g\geq 0$ or if $\Omega$ has a scalar curvature shield.  The same consideration was taken care of the fill-in problem as well.  

\begin{rem}
Let $3\leq n \leq 7$. From the perspective of non-compactness, if one has already assumed that on $H>0$ on $\S^{n-1}$, then one can allow $(\Omega, g)$ with $\p\Omega=\S$ to be an asymptotically flat manifold with non-negative scalar curvature to conclude that the Brown-York mass is non-negative.  The reason is that one can chop off the manifolds with a minimal surface.  It is because the positive mass theorem with corners holds true with the presence of minimal boundary.  Similarly, for $\MW$, if one has already assumed $H>|\pi(\nu,\cdot)|$ which implies $H>|tr_{\S}k|$. Then we can show $\MW$ to be non-negative for the case where $(\Omega, g,k)$ with $\p\Omega=\S$ is a spin asymptotically flat initial data set satisfying the dominant energy condition. The reason is one can chop off by MOTS/ MITS while the Dirac operator approcah of spacetime positive mass theorem with corners on spin manifolds allows MOTS/ MITS boundary (\cite{GHHP}).  
\end{rem}

As remarked by \cite{CLZ} and as what we have shown in the proof of Theorem \ref{RelaxedWY}, the Dirac operator approach is easily generalised to the case where singularities along a hypersurface exist.  

Now, we can state the following version of the collection of results in \cite{T2} taking completeness and dominant energy shields into account using the spacetime positive mass theorem with completeness or dominant energy shields shown in \cite{CLZ} as an ingredient in their corresponding proofs of \cite{T2}.  

\begin{thm}(cf. \cite{ST} Theorem 4.1)
For $n\geq 3$,  let $(\Omega^n, g, k)$ be a compact or complete spin initial data set with dominant energy condition or a spin initial data set with a dominant energy shield. Assume that $\S:=\SS$ has finitely many components. Let $H$ denote the mean curvature of $\S$ with respect to the outward normal $\nu$.  Suppose $\S$ can be isometrically embedded into $\R^{n}$ as a strictly convex closed hypersurface.  Denote the mean curvature of isometric embedding of $\S$ into $\R^n$ with respect to the outward normal by $H_0$. If $H > |\omega|$, where $\omega=\pi(\cdot,\nu)$, then 
$$\mathcal{W}(\S):=\int_\S H_0 - \left( H-|\omega| \right) \geq 0,$$ 
and equality implies that $\Sigma$ is connected and $\Omega$ is in Minkowski space.  
\end{thm}

\begin{thm} (cf. \cite{SWWZ} Theorem 1.3)
Let $D_{SB}:=(\mathbb{S}^{n-1}, \gamma, \alpha, H, \beta )$ be a spacetime Bartnik data set. If $\gamma$ is isotopic to $\gamma_{std}$ in $\M^q_{psc}(\mathbb{S}^{n-1}):=\{\eta: C^q \text{\,\,metrics on\,\,} \mathbb{S}^{n-1} {\,\,with\,\,} R_{\eta}>0 \}$, where $q\geq 5$, then there exists a constant $h_0=h_0(n,\gamma)>0$ such that 
if 
$$H-f>0\,\,\, \text{and\,\,\,\,}\, \int_{\mathbb{S}^{n-1}} H-f \, d\mu_{\gamma}>h_0,$$
where $f:=\sqrt{(\tr_{\S}\alpha)^2+|\beta|_{\gamma}^2}$,
then $D_{SB}$ cannot admit a fill-in which is a compact or complete spin initial data set with dominant energy condition or a spin initial data set with a dominant energy shield. 
\end{thm}

\begin{thm} (cf. \cite{SWWZ} Theorem 1.4)
Let $D_{SB}:=(\mathbb{S}^{n-1}, \gamma, \alpha, H, \beta )$ be a spacetime Bartnik data set. If $\gamma\in\M^n_{c,d}:=\{\eta: C^{\infty} \text{\,\,metrics on\,\,} \mathbb{S}^{n-1} {\,\,with\,\,}\\ |Rm_{\eta}|\leq c, \, diam({\eta})\leq d, \, vol(\eta)=vol(\gamma_{std})  \}$, then there exists a constant $C_0(n,c,d)>0$ such that if 
$$H-f\geq C_0,$$
where $f:=\sqrt{(\tr_{\S}\alpha)^2+|\beta|_{\gamma}^2}$,
then $D_{SB}$ cannot admit a fill-in which is a compact or complete spin initial data set with dominant energy condition or a spin initial data set with a dominant energy shield. 
\end{thm}

\begin{thm} (cf. \cite{SWW} Theorem 1.2)
Let $D_{SB}:=(\S^{n-1}, \gamma, \alpha, H, \beta )$ be a spacetime Bartnik data set where $\S^{n-1}$ can be smoothly embedded into $\R^n$ and $\gamma$ is smooth. There exists a constant $C_0=C_0(\Sigma,\gamma)>0$ such that if $$H-f\geq C_0,$$ 
where $f:=\sqrt{(\tr_{\S}\alpha)^2+|\beta|_{\gamma}^2}$,
then $D_{SB}$ cannot admit a fill-in which is a compact or complete spin initial data set with dominant energy condition or a spin initial data set with a dominant energy shield. 
\end{thm}

\end{document}